\documentclass[a4paper,10pt]{article}
\usepackage{amsmath,amscd,amsthm,amsfonts,amssymb,epsf,latexsym,wasysym}
\usepackage[all]{xy}
\usepackage{graphicx}

\makeatletter

\long\def\@makefnt#1{\parindent 1em\noindent
             \hb@xt@1.8em{\hss\@textsuperscript{}}#1}
\long\def\@ftntext#1{\insert\footins{%
     \reset@font\footnotesize
     \interlinepenalty\interfootnotelinepenalty
     \splittopskip\footnotesep
     \splitmaxdepth \dp\strutbox \floatingpenalty \@MM
     \hsize\columnwidth \@parboxrestore
     \color@begingroup
       \@makefnt{%
         \rule\z@\footnotesep\ignorespaces#1\@finalstrut\strutbox}%
     \color@endgroup}}%
\def\subjclass#1{%
   \@ftntext{2010 {\itshape Mathematics Subject Classification.}\enspace
#1.}}
\def\keywords#1{%
   \@ftntext{{\itshape Keywords.}\enspace #1.}}
\makeatother

\def\moins{\raise 1pt\hbox{{$\scriptstyle -$}}}
\def\plus{\raise 1pt\hbox{{$\scriptstyle +$}} }

\newtheorem{theorem}{Theorem}

\newtheorem{lemma}[theorem]{Lemma}

\newtheorem{observation}[theorem]{Observation}

\def\dim{\mathop{\rm dim}}

\def\C{{\bf C}}

\def\P{{\bf P}}
\def\Q{{\bf Q}}

\def\\Q{\widetilde{Q}}

\begin{document}

\title{\bf On some properties of the \L ojasiewicz exponent}
\author{Christophe Eyral, Tadeusz Mostowski and Piotr Pragacz}

\date{}

\maketitle

\begin{abstract}
In this note, we investigate the behaviour of the \L ojasiewicz exponent under hyperplane sections and its relation to the order of tangency.
\end{abstract}

\vskip 1cm

\noindent
\textbf{1. Introduction and statements of the results}\vskip 2mm

It is well known (see \cite{L,L1}) that any pair of closed analytic subsets $X, Y\subset \mathbf{C}^m$ ($m\geq 2$) satisfies so-called \emph{\L ojasiewicz regular separation property} at any point of $X\cap Y$. Precisely, for any $x^0\in X\cap Y$ there are constants $c,\nu >0$ such that for some neighbourhood $U\subset \mathbf{C}^m$ of $x^0$ we have
\begin{equation}\label{cond-L1}
\rho(x,X)+\rho(x,Y)\ge c \, \rho(x, X\cap Y)^{\nu} \quad \mbox{for $x\in U$},
\end{equation}
where $\rho$ is the distance induced by the standard Hermitian norm on $\mathbf{C}^m$.
Note that if $x^0\notin\mbox{int}(X\cap Y)$, where the interior is computed in $\mathbf{C}^m$, then necessarily $\nu\geq 1$ (see \cite{D}). Also, observe that $X$ and $Y$ satisfy \eqref{cond-L1} with a constant $\nu\geq 1$ if and only if there exist a neighbourhood $U'$ of $x^0$ and a constant $c'>0$ such that
\begin{equation}\label{cond-L2}
\rho(x,Y)\ge c' \rho(x, X\cap Y)^{\nu} \quad \mbox{for $x\in U'\cap X$}
\end{equation}
(see \cite{L,CT,D}).
Any exponent $\nu$ satisfying the relation \eqref{cond-L1} for some $U$ and $c>0$ is called a \emph{regular separation exponent} of $X$ and $Y$ at $x^0$. The infimum of such exponents is called the \emph{\L ojasiewicz exponent} of $X$ and $Y$ at $x^0$ and is denoted by ${\mathcal L}(X,Y;x^0)$; it is important to observe that the latter is a regular separation exponent itself (see \cite{S}). The number ${\mathcal L}(X,Y;x^0)$ is an interesting metric invariant of the pointed pair $(X,Y;x^0)$ which have been the subject of vast studies in analytic geometry (see, for instance, the references in \cite{S}). 

The goal of this note is to investigate the behaviour of the \L ojasiewicz exponent under hyperplane sections. Precisely we show the following theorem.

\begin{theorem}\label{ml}
Let $X$ and $Y$ be closed analytic subsets in $\mathbf{C}^m$, and let $x^0\in X\cap Y$ such that ${\mathcal L}(X,Y;x^0)\ge 1$.
Then for a general hyperplane $H_0$ of $\C^m$ passing through $x^0$ we have
\begin{equation*}
{\mathcal L}(X\cap H_0,Y\cap H_0;x^0)\le {\mathcal L}(X,Y;x^0).
\end{equation*}
\end{theorem}

This theorem is a consequence of the following result, which is the main part of the present work. 

\begin{theorem}\label{mc}
Let $X$ be a closed analytic subset in $\C^m$, and let $x^0\in X$. Then for a general hyperplane $H_0$ of $\C^m$ passing through $x^0$, there exist a constant $c>0$ and a neighbourhood $U$ of $x^0$ such that for all $x\in U\cap H_0$ we have
\begin{equation*}
\rho(x,X\cap H_0)\le c\,\rho(x,X).
\end{equation*}
\end{theorem}

Theorems \ref{ml} and \ref{mc} are proved in Sections 2 and 3 respectively. To conclude this paper, in Section 4, we also briefly discuss the relation between the \L ojasiewicz exponent and the order of tangency for pairs of closed analytic submanifolds of $\C^m$ with the same dimension.
\vskip 5mm

\noindent
\textbf{2. Proof of Theorem \ref{ml}} \vskip 2mm

Without loss of generality, we may assume that $x^0$ is the origin $0\in\C^m$.
If $\nu$ is a regular separation exponent for $X$ and $Y$ at $0$, then $\nu\ge \mathcal{L}(X,Y;0)\ge 1$, and by \eqref{cond-L2}, for some $c'>0$ we have
\begin{equation}\label{proof-fi}
\rho(x,Y)\ge c' \rho(x, X\cap Y)^{\nu}
\end{equation}
for all $x\in X$ near $0$. By Theorem \ref{mc}, applied to $X\cap Y$, for a general hyperplane $H_0$ of $\C^m$ there is a constant $c>0$ such that for all $x\in H_0$ near $0$ we have
\begin{equation*}
c\,\rho(x,X\cap Y)^{\nu}\ge \rho(x,X\cap Y\cap H_0)^{\nu}.
\end{equation*}
Combined with \eqref{proof-fi}, this gives
\begin{align*}
\rho(x,Y\cap H_0) \ge \rho(x,Y)
\ge c'\,\rho(x,X\cap Y)^{\nu} \ge (c'/c)\,\rho(x,X\cap Y \cap H_0)^{\nu}
\end{align*}
for all $x\in X\cap H_0$ near $0$, so that $\nu$ is a regular separation exponent for $X \cap H_0$ and $Y \cap H_0$ at $0$. Applying this with $\nu=\mathcal L(X,Y;x^0)$ shows that
\begin{equation*}
{\mathcal L}(X\cap H_0,Y\cap H_0;x^0)\le {\mathcal L}(X,Y;x^0).
\end{equation*}
\vskip 4mm

\noindent
\textbf{3. Proof of Theorem \ref{mc}}\vskip 2mm
 
It strongly relies on the Lipschitz equisingularity theory of complex analytic sets developed in \cite{M} by the second named author. Throughout, we always work with Hermitian orthonormal bases $\{e_1,\ldots,e_m\}$ in $\C^m$, and the corresponding coordinates $x=(x_1,\ldots,x_m)$.
As in Section 2, we assume $x^0=0$ and we work in a small neighbourhood of it.

Let $\check\P^{m-1}$ denote the set of all hyperplanes of $\C^m$ through $0$, with its usual structure of manifold. The distance between two elements $H,K\in\check\P^{m-1}$  is the angle $\varangle(H,K)$ between them, that is, 
\begin{equation*}
\varangle(H,K):=\arccos \frac{|\langle v,w \rangle|}{|v|\, |w|} \in [0,\pi/2]
\end{equation*}
where $v$ and $w$ are normal vectors to the hyperplanes $H$ and $K$, respectively, and $\langle \cdot,\cdot\rangle$ is the standard Hermitian product on $\C^{m}$ (see, e.g., \cite{Sc}). 
\vskip 4mm

\noindent
\emph{Step 1.} Let
$$
\mathcal{X}:=\{(H,x)\in \check\P^{m-1}\times \C^m\mid x\in H\cap X\}.
$$
By Proposition 1.1 of \cite{M}, in a neighbourhood 
$$
\mathcal U:=\{(H,x)\in \check\P^{m-1}\times \C^m\mid \varangle(H_0,H)<a \mbox{ and } \vert x\vert<b\}
$$
of a generic $(H_0,0)$, we have that $\mathcal{X}$ is \emph{Lipschitz equisingular} over $\check\P^{m-1}\times \{0\}$. That is, for any  $(H,0)\in {\mathcal U} \cap (\check\P^{m-1}\times \{0\})$, there is a (germ of) Lipschitz homeomorphism 
\begin{equation*}
\varphi\colon (\check\P^{m-1}\times \C^m, (H,0)) \to (\check\P^{m-1}\times \C^m, (H,0))
\end{equation*}
(with a Lipschitz inverse)
such that $p\circ\varphi=p$ and $\varphi(\mathcal{X})=\check\P^{m-1} \times(H\cap X)$ (as germs at $(H,0)$). (Here, $p\colon \check\P^{m-1}\times \C^m\to \check\P^{m-1}$ is the standard projection.)
Actually, if $h=(h_1,\ldots,h_{m-1})$ are coordinates in $\check\P^{m-1}$ around $H_0$ such that 
$$
h_1(H_0)=\cdots =h_{m-1}(H_0)=0\,,
$$
then, locally near $(H_0,0)$, the standard ``constant'' vector fields $\partial_{h_j}$ ($1\leq j\leq m-1$) on $\check{\mathbf{P}}^{m-1}\times \{0\}$ can be lifted to Lipschitz vector fields $v_j$ on $\check{\mathbf{P}}^{m-1}\times \mathbf{C}^m$ such that the flows of $v_j$ preserve $\mathcal{X}$ (see the proof of Proposition 1.1 of \cite{M}, p.10).
So, in particular, $v_j$ is a Lipschitz vector field of the form
\begin{equation*}
v_j(h,x)=\partial_{h_j}(h,x)+\sum_{\ell=1}^m w_{j\ell}(h,x) \,
\partial_{x_\ell}(h,x),
\end{equation*}
so that $v_j(h,0)=\partial_{h_j}(h,0)$ and there exists a constant $c'>0$ such that
\begin{equation}\label{proof-mc-ineg}
\vert w_{j\ell}(h,x)\vert\le c'\,\vert x\vert \mbox{ near $0$}
\end{equation}
for all $j,\ell$. 
\vskip 4mm

\noindent
\emph{Step 2.} Pick a point $y^0\in H_0$. We want to prove that if $y^0$ is sufficiently close to $0$, then 
\begin{equation}\label{iop}
\rho(y^0,X\cap H_0)\le c\,\rho(y^0,X)
\end{equation}
 for some constant $c>0$ independent of $y^0$.
Let $y^1\in X$  be one of the closest points to $y^0$, that is, $\rho(y^0,X)=\vert y^1-y^0\vert$.
If $y^0\in X$, then $\rho(y^0,X\cap H_0)=\rho(y^0,X)=0$, and the inequality \eqref{iop} is obviously true. So, hereafter, we assume that $y^0\notin X$.
Of course, without loss of generality, we may also assume that $|y^0|<b$ and $|y^1|< b$.
Choose $H_1\in \check\P^{m-1}$ such that $y^1\in H_1$ and $\varangle(H_0,H_1)$ is minimal. 
If $\varangle(H_0,H_1)=0$ (i.e., if $y^1\in H_0$), then again $\rho(y^0,X\cap H_0)=\rho(y^0,X)$ and \eqref{iop} is true. From now on, let us assume that $\varangle(H_0,H_1)\not=0$.
Then we have the following lemma.

\begin{lemma}\label{aa}
If $(H_1,y^1)\notin \mathcal U$ (i.e., if $\varangle (H_0,H_1)\geq a$), then there exists $a'>0$ depending only on $a$ such that
\begin{equation*}
\vert y^1-y^0\vert \ge a'\, \vert y^0\vert.
\end{equation*}
\end{lemma}

In particular, since $0\in X\cap H_0$, if $(H_1,y^1)\notin \mathcal U$ then we have
\begin{equation}\label{inegsou}
\rho(y^0,X\cap H_0)\le \vert y^0\vert\le (1/a')\, \rho(y^0,X)
\end{equation}
as desired.
\vskip 2mm

\noindent
\emph{Proof of Lemma \ref{aa}.} 
By a proper choice of the basis $\{e_1,\ldots,e_m\}$, we may assume that $H_0$ is defined by the equation $x_m=0$, so that $e_m$ is orthogonal to $H_0$. Now, if $x_m=\sum_{\ell=1}^{m-1} q_\ell \, x_\ell$ is an equation for $H_1$, then, clearly, for each $1\leq \ell\leq m-1$, the vector $E_\ell:=e_\ell+q_\ell e_m$ is in $H_1$. Thus, if  $N=\sum_{\ell=1}^{m-1}u_\ell e_\ell+u_m e_m$ is a normal vector to $H_1$, then we must have $\langle N, E_\ell\rangle=0$, and hence, $u_\ell=-u_m \bar q_\ell$, so that we can take $N:=-\sum_{\ell=1}^{m-1}\bar q_\ell\, e_\ell+e_m$.

Now, saying that $\varangle(H_0,H_1)$ is minimal means that 
\begin{equation*}
\cos\varangle(H_0,H_1)=\frac{|\langle N, e_m\rangle|}{|N|\, |e_m|}=\frac{1}{\sqrt{1+\sum_{\ell=1}^{m-1} |q_\ell|^2}}
\end{equation*}
 is maximal, that is, $\sum_{\ell=1}^{m-1} |q_\ell|^2$ is minimal.
By adjusting the choice of the basis, we may further assume that $y^1=(y^1_1,0,\ldots,0,y^1_m)$, so that its orthogonal projection onto $H_0$ is $y^2:=(y^1_1,0,\ldots,0)$. As $y^1\in H_1$, we have $q_1=y^1_m/y^1_1\not=0$. Thus, $\sum_{\ell=1}^{m-1} |q_\ell|^2$ is minimal if and only if $q_2=\cdots=q_{m-1}=0$. So, if $\varangle(H_0,H_1)$ is minimal, then $H_1$ is given by the equation $x_m=q_1 x_1$.

\begin{figure}[t]
\centerline{\includegraphics[width=6cm,height=2.4cm]{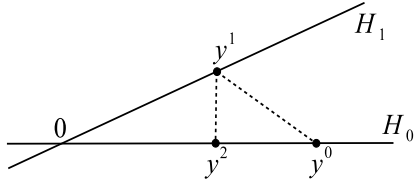}}
\caption{Hyperplanes $H_0$ and $H_1$}
\label{Fig1}
\end{figure}

It follows that if $\varangle (H_0,H_1)\geq a$ (assumption of the lemma), then we must have 
\begin{equation*}
\cos\varangle(H_0,H_1)=1/\sqrt{1+|q_1|^2}\leq a_1,
\end{equation*}
and hence $|q_1|\geq a_2$,
for some constants $a_1,a_2>0$ depending only on $a$.
Now, clearly, we may always assume $|y^0-y^1|< (1/10)\, |y^0|$. Thus, $|y^2-y^0|\leq |y^1-y^0|<(1/10)\, |y^0| $, and hence,
\begin{equation*}
|y^2-0|=|y^1_1|>(9/10)\, |y^0|
\end{equation*}
(see Figure \ref{Fig1}).
It follows that
\begin{equation*}
|y^0-y^1|\geq |y^1-y^2|=|q_1|\, |y^1_1|\geq a_2\, (9/10)\, |y^0|,
\end{equation*}
and this completes the proof of Lemma \ref{aa}.

\vskip 4mm

\noindent
\emph{Step 3.} 
Lemma \ref{aa} solves the case where $(H_1,y^1)\notin \mathcal U$ (see \eqref{inegsou}).
Now let us look at the case where $(H_1,y^1)\in \mathcal U$; here comes Lipschitz equisingularity (see Step 1). Let $h^1=(h_1^{1},\ldots,h_{m-1}^{1})$ be the coordinates of $H_1$. (Note that $|h^1|\leq d \cdot\varangle(H_0,H_1)$ for some constant $d>0$ independent of $H_1$.)
Consider the Lipschitz vector field $v$ on $\check\P^{m-1}\times \C^m$ defined by
\begin{align*}
v(h,x) & :=-\sum_{j=1}^{m-1}h_j^1\, v_j(h,x) \\
& = -\sum_{j=1}^{m-1}h_j^1\, \partial_{h_j}(h,x) +
\sum_{\ell=1}^m \bigg( -\sum_{j=1}^{m-1}h_j^1\,  w_{j,\ell}(h,x) \bigg) \, \partial_{x_\ell}(h,x),
\end{align*}
and look at the integral curve $\gamma(t)=(h(t),x(t))$ 
of  $v$ starting at $(H_1,y^1)$. So, in particular, we have:
\begin{align*}
& \dot h_j(t)=-h^1_j, \quad \dot x_\ell(t)=-\sum_{j=1}^{m-1}h_j^1\,  w_{j,\ell}(h,x),\\
& h_j(0) =h^1_j, \quad \ \, x_\ell(0)=y^{1}_{\ell}.
\end{align*}
As the flows of the vector fields $v_j$ preserve $\mathcal{X}$ and since $\gamma(0)\in\mathcal{X}$, the curve $\gamma(t)$ lies in  $\mathcal{X}$. Moreover, since $h_j(t)=h^1_j(1-t)$, we have $h_j(1)=0$ for all $j$, and hence $x(1)$ lies in $H_0$. 
Finally, observe that the length $L_I(x)$ of the restriction of the curve $x(t)$ to the compact interval $I=[0,1]$ satisfies
\begin{align*}
L_I(x):=\int_0^1 \vert \dot x(t)\vert  \, dt & \le
c_1\int_0^1 \sum_{j=1}^{m-1}\bigg(\vert h^1_j\vert\cdot \bigg(\sum_{\ell=1}^m 
\vert w_{j,\ell}(\gamma(t))\vert\bigg)\bigg) dt \\ 
& \overset{\mbox{\tiny by \eqref{proof-mc-ineg}}}{\le} c_2\, \vert h^1\vert \int_0^1\vert x(t)\vert\, dt
\le c_3\, \vert h^1\vert\, \vert x(0)\vert\le c_4\, \vert y^0-x(0)\vert
\end{align*} 
for some constants $c_i>0$ independent of $y^0$, $H_1$ and $y^1$. 
The first and third inequalities are clear. The second one follows from the crucial relation \eqref{proof-mc-ineg} (i.e., from Lipschitz equisingularity).
To show the last inequality, we may proceed as in the proof of Lemma \ref{aa}, exchanging the roles of $H_0$ and $H_1$. Namely, for a new proper choice of the basis, we may assume that $H_1$ is defined by $x_m=0$ and that $y^0=(y^0_1,0,\ldots,0,y^0_m)$, so that the orthogonal projection of $y^0$ onto $H_1$ is $y^3:=(y^0_1,0,\ldots,0)$. As the angle $\varangle(H_0,H_1)$ is minimal, we may suppose that $H_0$ is given by an equation of the form $x_m=q_1 x_1$. Clearly, we may also assume that $|y^0-y^1|< (1/10)\, |y^1|$. Thus $|y^3-y^1|\leq |y^0-y^1|< (1/10)\, |y^1|$, and hence, $|y^3-0|=|y^0_1|> (9/10)\, |y^1|$. It follows that 
\begin{align*}
|y^0-y^1|\geq |y^0-y^3|=|q_1|\, |y^0_1|>|q_1|\, (9/10)\, |y^1|.
\end{align*} 
But we have
\begin{align*}
|h^1|\leq d \cdot \varangle(H_0,H_1)=d \cdot \arccos (1/\sqrt{1+|q_1|^2})\leq d'\, |q_1|,
\end{align*} 
where the constants $d,d'>0$ are independent of $y^0$, $H_1$ and $y^1$. It follows that 
\begin{align*}
|y^0-y^1|\geq \frac{9}{10\, d'}\, |h^1|\, |y^1|
\end{align*} 
 as desired (remind that $y^1=x(0)$). Now, by the estimate of the length $L_I(x)$ given above, we have
\begin{align*}
\rho(y^0,X\cap H_0) & \le \vert y^0-x(1)\vert\le \vert y^0-x(0)\vert +\vert x(0)-x(1)\vert \le \vert y^0-x(0)\vert + L_I(x)\\ 
& \le (1+c_4)\, \vert y^0-x(0)\vert = (1+c_4)\, \rho(y^0,X),
\end{align*}
and this completes the proof of Theorem \ref{mc}.

\vskip 2mm

\noindent
\emph{Remark.} 
Note that the proof of Theorem \ref{mc} (and hence of Theorem \ref{ml}) given above only depends on the Lipschitz equisingularity theory of complex analytic sets developed in \cite{M} by the second named author. Real versions of this theory for the semi-analytic and subanalytic categories were addressed by A.~Parusi\'nski in \cite{P4,P1,P2,P3} while the case of sets definable in a polynomially bounded o-minimal structure was obtained by Nguyen Nhan and G. Valette in \cite{NV}. Theorems \ref{ml} and \ref{mc} must then be true in these categories as well.

\vskip 5mm

\noindent
\textbf{4. Remark on the \L ojasiewicz exponent and the order of tangency}\vskip 2mm

To conclude this paper, we give a lower bound for the \L ojasie\-wicz exponent $\mathcal L(X,Y;x^0)$ of two $p$-dimensional closed analytic submanifolds $X$ and $Y$ of $\C^m$ at $x^0\in X\cap Y$ in terms of the order of tangency of $X$ and $Y$ at $x^0$. 

Following \cite{J,DMP}, we say that the order of tangency between $X$ and $Y$ at $x^0$ is greater than or equal to an integer $k$ if there exist parametrizations (i.e., biholomorphisms onto their images)
$$
q\colon (U,u^0) \to (X,x^0) 
\quad\mbox{and}\quad
q'\colon (U,u^0) \to (Y,x^0),
$$
where $U \ni u^0$ is an open subset of \,$\C^p$, such that
\begin{equation}\label{oog}
q(u)-q'(u)=o(\vert u-u^0\vert^k)
\end{equation}
when $U\ni u\to u^0$. The \emph{order of tangency} between $X$ and $Y$ at $x^0$ (denoted by $s(X,Y;x^0)$) is the supremum of all such integers $k$.  

\begin{observation}\label{sl} 
Let $X$ and $Y$ be $p$-dimensional closed analytic submanifolds of $\C^m$, and let $x^0\in X\cap Y$. Suppose that $s(X,Y;x^0)$ is finite. If $\mathcal{L}(X,Y;x^0)\geq 1$, then 
\begin{equation*}
s(X,Y;x^0)\leq \mathcal{L}(X,Y;x^0)-1.
\end{equation*}
\end{observation}

\noindent
\emph{Proof.}
Put $s:=s(X,Y;x^0)$, $\mathcal{L}:=\mathcal{L}(X,Y;x^0)$, and for this proof write
$\mathbf{C}^m=\mathbf{C}_x^p\times \mathbf{C}_y^{m-p}$ where $x=(x_1,\ldots,x_p)$ and $y=(x_{p+1},\ldots,x_m)$. As above, we assume that $x^0$ is the origin $0\in \C^m$.
In a neighbourhood of $0$, the analytic submanifold $X$ is given by $y=f(x)$ for some analytic function 
\begin{equation*}
f=(f_1,\ldots,f_{m-p})\colon\, (\mathbf{C}_x^p,0)\to(\mathbf{C}_y^{m-p},0).
\end{equation*}
 Similarly, $Y$ is also the graph of an analytic function $g$, and without loss of generality, we may assume that $g=0$. Now, let $s'$ be the smallest integer $k$ for which there exists a multi-index $\alpha=(\alpha_1,\ldots,\alpha_p)$ such that $|\alpha|=\alpha_1+\cdots+\alpha_p=k$ and $D^{\alpha}(f-g)(0)\not=0$. Clearly, $s=s'-1$. 
Each component $f_i$ has the Taylor expansion
$$
f_i(x)=F_i(x)+o(|x|^{r_i})
$$
where $F_i$ is a homogeneous polynomial of degree $r_i$. Of course, we may assume $r_1\leq r_i$ for all $i$, so that $r_1=s'$. 
Consider the standard projection 
\begin{align*}
\pi\colon\, \mathbf{C}_x^p\times \mathbf{C}_y^{m-p}\to \mathbf{C}_x^p,
\end{align*}
 and look at the hypersurface $\pi(X\cap Y)=\{x\in \mathbf{C}_x^p\, ;\, f(x)=0\}$ of $\mathbf{C}_x^p$. It is easy to see that if $L$ is a line through~$0$ which is not in the tangent cone of $\pi(X\cap Y)$ at $0$, then
\begin{align*}
\rho(x,\pi(X\cap Y)) \sim |x|
\end{align*}
for $x\in L$ near $0$.\footnote{As usual, the expression $\varphi(x)\sim\psi(x)$ for $x\in E$ near $0$ means that there exist constants $c,c'>0$ such that $c\, \psi(x)\leq \varphi(x)\leq c'\, \psi(x)$ for all $x\in E$ near $0$.} 
Now, if $F_1\not=0$ on $L$, then for any $x\in L$ near $0$, we also have
\begin{align*}
|f_1(x)|\sim |x|^{r_1}=|x|^{s'} \quad \mbox{and} \quad
|f_i(x)|\leq a\, |x|^{r_i} \leq a\, |x|^{s'}
\end{align*}
for some constant $a>0$.
It follows that for any $(x,y)\in \pi^{-1}(L)\cap X=\{(x,y)\, ;\, x\in L\mbox{ and } y=f(x)\}$ near~$0$, we have
\begin{align*}
\rho((x,y),Y)=|f(x)|\sim |x|^{s'} \quad\mbox{and}\quad
\rho((x,y),X\cap Y) \sim |x|.
\end{align*}
Now, the \L ojasiewicz exponent $\mathcal{L}$ satisfies
$\rho((x,y),Y)\geq c\, \rho((x,y),X\cap Y)^{\mathcal{L}}$, that is,
$|x|^{s'}\geq c\, |x|^{\mathcal{L}}$ for some constant $c>0$. Thus $s'\leq \mathcal{L}$, and hence, $s=s'-1\leq\mathcal{L}-1$.

\vskip 2mm
\noindent\emph{Remark.} We may also investigate the relationship between $s:=s(X,Y;x^0)$ and $\mathcal{L}:=\mathcal{L}(X,Y;x^0)$ using Theorem \ref{ml} but this second approach only gives the inequality $s<\mathcal{L}$. However, for completeness, let us briefly explain the argument.
First, we consider the special case where $x^{0}$ is an isolated point of $X\cap Y$. In this case, there exists a constant $c'>0$ such that
\begin{equation*}
\rho(x,Y) \geq c'\, \rho(x,X\cap Y)^{\mathcal{L}}=c'\, |x-x^0|^{\mathcal{L}}
\quad\mbox{for } x\in X \mbox{ near } x^0,
\end{equation*}
or equivalently, $\rho(q(u),Y) \geq c'\, |q(u)-q(u^0)|^{\mathcal{L}}$ for $u$ near $u^0$.
Since $q$ is locally bi-Lipschitz, there exists a constant $c''>0$ such that
\begin{equation*}
c'\, |q(u)-q(u^0)|^{\mathcal{L}}\geq c''\, |u-u^0|^{\mathcal{L}}
\quad\mbox{for } u \mbox{ near } u^0.
\end{equation*}
Now, by \eqref{oog}, we have
\begin{equation*}
\rho(q(u),Y) \leq |q(u)-q'(u)| <c''\, |u-u^0|^s
\quad\mbox{for } u \mbox{ near } u^0.
\end{equation*}
Combining these relations gives 
\begin{equation*}
c''\, |u-u^0|^{\mathcal{L}}\leq \rho(q(u),Y)< c''\, |u-u^0|^s
\quad\mbox{for } u \mbox{ near } u^0,
\end{equation*}
and hence $s<\mathcal{L}$.

The general case (i.e., $\dim X\cap Y=n>0$) follows from the $0$-dimensional case and Theorem \ref{ml}. Indeed, take $n$ general hyperplanes $H_1,\ldots,H_n$ in $\C^m$ passing through $x^0$, so that $X\cap Y\cap H_1\cap \cdots \cap H_n$ is an isolated intersection. Let $s_i$ (respectively, ${\mathcal L}_i$) denote the order of tangency (respectively, the \L ojasiewicz exponent) of $X\cap H_1\cap \cdots \cap H_i$ and $Y\cap H_1\cap \cdots \cap H_i$ at $x^0$. Clearly, \eqref{oog} implies $s_i\leq s_{i+1}$ while Theorem \ref{ml} shows ${\mathcal L}_i\ge {\mathcal L}_{i+1}$. (Note that since $\mbox{int}(X\cap Y\cap H_1\cap \cdots \cap H_i)=\emptyset$, we have ${\mathcal L}_i\geq 1$, so that Theorem \ref{ml} applies.) Now the relation $s<\mathcal{L}$ follows from the inequality $s_n< {\mathcal L}_n$ ($0$-dimensional case).

\vskip 4mm

{\it Acknowledgments.} We warmly thank Tadeusz Krasi\'nski and the referee for valuable comments and suggestions which enabled us to improve the paper.

{\small \obeylines \parindent0pt
{\ }
Christophe Eyral
Institute of Mathematics, Polish Academy of Sciences
\'Sniadeckich 8, 00-656 Warszawa, Poland
E-mail: cheyral@impan.pl
{\ }
Tadeusz Mostowski
Department of Mathematics, Warsaw University
Banacha 2, 02-097 Warszawa, Poland
E-mail: tmostows@mimuw.edu.pl
{\ }
Piotr Pragacz
Institute of Mathematics, Polish Academy of Sciences
\'Sniadeckich 8, 00-656 Warszawa, Poland
E-mail: P.Pragacz@impan.pl
}

\end{document}